\def\R{{\rm I}\!{\rm R}}

\def\carre{\hbox{\vrule \vbox to 7pt{\hrule width 6pt \vfill \hrule}\vrule }}

\font\letitre = cmbx10 at 15pt

\font\leresume = cmbx10 at 10pt

\font\lessections = cmbx10 at 12pt

\font\lesauteurs = cmr10 at 11pt

\parindent = 0cm

\centerline {\letitre  Weyl composition of symbols in large dimension}

\bigskip
\bigskip

\centerline { \lesauteurs L. Amour and J. Nourrigat}

\bigskip

\centerline {\it Universit\'e de Reims, France}

\bigskip
\bigskip

\centerline{\leresume  Abstract}

\bigskip

This paper is concerned with the Weyl composition of symbols in large dimension.
We specify a class of symbols in order to estimate the Weyl symbol
of the product of two Weyl $h-$pseudodifferential operators,
with constants independent of the dimension. The proof includes a regularized
and a hybrid compositions together with a decomposition formula. We also analyze
in this context the remainder term of the semiclassical expansion of the Weyl composition.

\bigskip

\bigskip


\smallskip


\vskip 1cm

{\lessections 1. Statement of results.}

\bigskip

Composition of Weyl pseudodifferential operators is a largely
studied area in the literature, giving numerous classical results,
depending on the class where the symbols of the operators belong.

\bigskip

The purpose here is to establish estimates for the Weyl composition
of symbols, independently of the dimension, allowing in particular
the dimension to go to infinity. To this aim,  the two composed
symbols are chosen in a simple class, defined in such a way that the
constants appearing in the inequalities are also independent of the
dimension.

\bigskip

In a recent work with L. Jager [A-J-N-2], we obtain an upper bound
in the $L^2$ norm of operators with a symbol belonging to this
class. The constants involved in the inequality estimating this norm
are also independent of the dimension.

\bigskip

When  $A$ and $B$ are two functions on  $\R^{2n}$, bounded together with their derivatives, their Weyl composition, depending on the parameter $h>0$ (the Weyl symbol of the product of the two Weyl $h-$pseudodifferential operators with symbols  $A$
and $B$) is the function $C_h(A , B)$ formally defined on
$\R^{2n}$ by:
$$ C_h ^{weyl} (A , B) (X) = ( \pi h)^{-2n} \int _{\R^{4n}} A(X+ Y) B(X+Z) e^{-{2i \over h}
\sigma ( Y , Z) } dY dZ\leqno (1.1)$$
where $\sigma $ is the symplectic form ($\sigma (X , Y) =  y\cdot
\xi - x\cdot \eta$ for $X = (x , \xi)$ and $Y= (y , \eta)$ in
$\R^{2n}$). The theory may be found in H\"ormander [HO] Chapter 18,
(also see [LER], [S] and, in the semiclassical setting, see e.g.
[M], [R]). If $A$ and $B$ are bounded continuous functions on
$\R^{2n}$, one notes that equality (1.1) formally defines a tempered
distribution on $\R^{2n}$. If $A$ and $B$ are in the class $C^m$
($m$ being sufficiently large) with bounded derivatives up to order
$m$ then $C_h^{weyl}(A , B)$ is a bounded continuous function.

\bigskip

The objective of this work is to derive
estimates for  $ C_h ^{weyl}
(A , B)$ where all the constants are independent of the dimension $n$.
In order to do that, we shall first define a class of symbols where all
the constants are also specified.

\bigskip

{\bf Definition 1.1.} {\it Let $(\rho_j)_{(j\geq 1)}$ and
$(\delta_j)_{(j\geq 1)}$ be two sequences of real numbers $\geq 0$. Fix
$M\geq 0$ and an integer  $m\geq 0$. Define $S_m(M, \rho , \delta
)$ as the set of  functions $F $ continuous on  $\R^{2n}$ ($n\geq
1$) such that, for each multi-index $(\alpha , \beta)$ verifying
$\alpha _j\leq m$ and $\beta _j\leq m$ for all $j \leq n$, the derivative $\partial_x^{\alpha }
\partial_{\xi}^{\beta } F$ exists, is continuous and bounded, and satisfies:
$$\sup _{X\in \R^{2n}}  \left | \partial_x^{\alpha } \partial_{\xi}^{\beta } F(X)\right |\leq M
\prod_{j\leq n}  \rho _j^{\alpha_j}  \delta _j^{\beta_j}  \leqno
(1.2)$$

}

\bigskip

In [A-J-N-2] we give a precise upper bound of the $L^2$ norm of Weyl
$h-$pseudodifferential operators associated with a symbol $A$ in
$S_2 (M, \rho , \delta)$ when $h \rho_j \delta_j \leq 1$ for all
$j\leq n$. In the works of Bernard Lascar [LA-1] to [LA-4] one finds
an extensive analysis of pseudodifferential operators in large and
infinite dimension.

\bigskip

{\bf Theorem 1.2.} {\it There exists a universal constant
$K>0$ such that, for all $n\geq 1$, for any $A$ in $S_6(M,
\rho , \delta)$ and $B$ in $S_6(M', \rho , \delta )$, the Weyl composition $C_h(A , B)$
is a bounded function on $\R^{2n}$ and satisfies, if
$h \rho_j \delta_j \leq 1$ for all $j\leq n$:
$$ \sup_{X\in \R^{2n}} |C_h(A ,  B)(X)|\leq M'' \hskip 2cm M'' = M M' \prod _{j\leq n}
(1 + K h \rho_j \delta_j) \leqno (1.3)$$
If $A$ is in $S_{m} (M, \rho , \delta )$ and $B$ in $S_{m} (M',
\rho , \delta )$ ($m\geq 6$) then $C_h(A ,  B)$ belongs to $S_{m-6}
(M'', 2 \rho , 2 \delta )$, with $M''$ defined in (1.3).

}

\bigskip

 Next, we give the asymptotic expansion of the Weyl composition with constants again independent of the dimension.

\bigskip

{\bf Theorem 1.3.}  {\it For every $N \geq 1$, let $R_N $ be the function defined by:
$$ C_h ^{weyl} (A , B) (X) =  \sum _{k=0}^{N-1}
{h^k \over (2i)^k k!} \sigma (\nabla _Y , \nabla _Z) ^k [A(X+Y)
B(X+Z)]\Bigg |_{Y= Z = 0} + R_N(X ,  h) \leqno (1.4)$$
Then we have, for all $A$  in $S_{m} (M, \rho , \delta )$ and $B$
in $S_{m} (M', \rho , \delta )$ ($m\geq N+6$):
$$ R _N ( \cdot , h) \in S_{m - N - 6}\left ( M''{h^N  \over N! }  \left [
\sum _{j=1}^n \rho_j \delta _j\right ]^N , 2 \rho , 2 \delta \right
) \leqno (1.5)$$
where $M''$ is defined in (1.3).

 }

\bigskip

The idea of the proof is to first introduce a regularized composition
 $ C_h^{reg} (A , B)$ for any functions $A$ and $B$
bounded on $\R^{2n}$. Namely, it is defined as the Wick symbol of the product of the two operators with anti-Wick symbols $A$ and $B$ respectively.
The $L^{\infty}(\R^{2n})$ norm of this regularized composition  is bounded by the product of the $L^{\infty}(\R^{2n})$ norms  of $A$ and $B$
(see Section 2). Immediately thereafter, we define for all subset  $I$
in $\{ 1 , ... , n\}$, a hybrid composition $ C_h^{hyb, I} (A , B)$ behaving as a Weyl composition with respect to the variables
$x_j$ with $j\in I$, and behaving as a regularized
composition with respect to the variables $x_j$ with $j\in
\{1 , ... , n\} \setminus I$.
In the next step, on the basis of a decomposition of the identity, Proposition 3.1 provides a decomposition of the Weyl composition $ C_h ^{weyl} (A , B)$ as a sum, where each term in the sum is related to a hybrid composed symbol associated with some subset $I$ of $\{ 1 , ... , n\}$, the sum being taken over all these subsets. As a further step, the hybrid composition is written in Proposition 2.1 as an integral expression,
 then, integrations by parts combined with other techniques allow to bound these hybrid compositions (Section 4).
In the last step, it  remains to take the sum of the bounds associated to each subset $I$ of $\{ 1 , ... , n\}$ to derive the estimate
(1.3). The other claim in Theorem 1.2 and Theorem 1.3 are then deduced relying on standard arguments (Section 5).

\bigskip

{\lessections 2. Regularized and hybrid compositions of symbols.}

\bigskip

We shall first study a composition law on $L^{\infty}
(\R^{2n})$ which shall be a continuous bilinear map. For that purpose, we define for all $A$ and $B$ in $L^{\infty } (\R^{2n})$, for all $X$ in $\R^{2n}$:
$$ C_h^{reg} (A , B) (X) =e^{{h\over 4}\Delta} C_h ^{weyl} ( e^{{h\over
4}\Delta}A , e^{{h\over 4}\Delta} B)(X) \leqno (2.1)$$
Then, for all subsets $I$ in $\{1, ... , n\}$, we also define:
$$ C_h^{hyb, I} (A , B) (X) =e^{{h\over 4}\Delta_{I^c}} C_h^{weyl}  ( e^{{h\over
4}\Delta_{I^c}}A , e^{{h\over 4}\Delta_{I^c}} B)(X)\leqno (2.2)$$
where $I^c$ is the complement of $I$ in $\{1, ... , n\}$ and
$$\Delta_{I^c} = \sum _{j\in I^c} {\partial ^2 \over \partial x_j^2}
+{\partial ^2 \over \partial \xi_j^2}\leqno (2.3)$$
Thus, if $I = \emptyset $ then $ C_h^{hyb, I} (A , B)=  C_h^{reg} (A ,
B)$ and if $I =\{1, ... , n\}$ then $ C_h^{hyb, I} (A , B)=  C_h^{weyl}
(A , B)$.

\bigskip

For all subsets $I$ of $\{1, ... , n\}$ and for every functions $A$ and $B$ on $(\R^2)^{I}$,   $C_h^{reg, I}(A ,
B)$ denotes the function on $(\R^2)^{I}$ defined as in  (2.1), when replacing
 $\{1 , \dots, n\}$ by $I$. For all functions $A$ on
$\R^{2n}$ and for all $X_I$ in $(\R^2)^{I}$, we define a function $A_{X_I } $ on
 $(\R^2)^{I^c}$ setting $A_{X_I } (X_{I^c}) = A (X_I , X_{I^c})$.
 With these notations, we may write:
 $$ C_h^{hyb, I} (A , B) (X_I , X_{I^c})=( \pi h)^{-2|I|} \int _{(\R^{4})^I}
 C_h ^{reg, I^c} (A_{X_I+ Y_I},  B_{X_I+ Z_I}) ( X_{I^c})   e^{-{2i \over h}
\sigma ( Y_I , Z_I) } dY_I dZ_I \leqno (2.4)$$

\bigskip

We shall now express $ C_h^{hyb, I} (A , B)$ under an integral form.

\bigskip

{\bf Proposition 2.1.} For each subset $I$  of $\{1 ,
..., n\}$ we have:
$$ C_h^{hyb, I} (A , B) (X) = \int_{\R^{4n}} A(X + Y) B(X+Z) K_{I, h} ( Y , Z)
d\lambda (Y, Z)  \leqno (2.5)$$
$$ K_{I, h} ( Y , Z) = (\pi h)^{-2|I|} (2 \pi h)^{-2|I^c|}
e^{-{2i \over h} \sigma ( Y_I , Z_I) } e^{{1\over 2h}Z_{I^c}\cdot
\overline {Y_{I^c}}} e^{-{1\over 2h} (|Y_{I^c}|^2 + |Z_{I^c}|^2)}
\leqno (2.6)$$

\bigskip

{\it Proof.}  We first prove the proposition
for $I= \emptyset$, that is to say, for the function $C_h^{reg }(A ,
B)$. Let us recall the coherent states:
$$ {\Psi}_{X , h} (u)  =  {\Psi}_{a , b , h} (u) =  (\pi h)^{
-{n/4}} e^{-{| u-a|^2 \over 2h}} e^{{i \over h} u .b - {i \over 2h}
a. b} \hskip 2cm u\in \R^n  \leqno (2.7) $$
We denote by $ Op_h^{AW} (A)$  the  anti-Wick operator associated with the
symbol  $A$, that is to say, the operator defined for all $f$ and
$g$ in $L^2(\R^n)$ by:
$$< Op_h^{AW} (A) f, g> = (2\pi h)^{-n} \int_{\R^{2n}} A(X) < f ,{\Psi}_{X ,
h}> < {\Psi}_{X , h} , g> dX\leqno (2.8)$$
We know that the Weyl symbol of this operator is $e^{{h\over
4}\Delta }A$. Consequently, $C_h ( e^{{h\over 4}\Delta}A ,
e^{{h\over 4}\Delta} B)$ is the Weyl symbol of the product $
Op_h^{AW} (A) \circ  Op_h^{AW} (B)$. Moreover, we call Wick symbol of
an operator $C$ bounded in $L^2(\R^n)$, the function
defined on $\R^{2n}$ by:
$$ \sigma_h^{wick}(C) (X) = < C \Psi_{Xh}, \Psi _{X , h}> \hskip 2cm
X\in \R^{2n}\leqno (2.9)$$
If $C$ is written under the form $C = Op_h^{weyl }(F)$, we know that its Wick symbol $\sigma_h^{wick}(C) =e^{{h\over 4}\Delta}
F$. These points imply that:
$$ C_h^{reg} (A , B)=  \sigma_h^{wick}( Op_h^{AW}(A) \circ
Op_h^{AW}(B))\leqno (2.10)$$
Taking these considerations into account, it appears:
$$C_h^{reg} (A , B) (X) = ( 2\pi h)^{-2n} \int_{\R^{4n}}
A(Z_1) B(Z_2) < \Psi_{X , h} ,\Psi_{Z_2 , h} >  < \Psi_{Z_2 , h}
,\Psi_{Z_1 , h} >  < \Psi_{Z_1 , h} ,\Psi_{X , h} > dZ_1dZ_2$$
Then recalling:
$$ <{\Psi}_{X , h} , {\Psi}_{Y , h} > = e^{-{1\over 4h}|X-Y|^2}
e^{{i\over 2h}\sigma (X , Y)}\leqno (2.11)$$
we express $C_h^{reg} (A , B)$ as:
$$C_h^{reg} (A , B) (X) = ( 2\pi h)^{-2n} \int_{\R^{4n}} A(X+Y_1)
B(X+Y_2) e^{{1\over 2h}Y_2\cdot \overline {Y_1} } e^{-{1\over 2h}
(|Y_1|^2 + |Y_2|^2)} dY_1dY_2$$
We similarly derive an analogous equality for $C_h^{reg ,
I^c}$ substituting $\{ 1 , ... , n\}$ to $I^c$. We then deduce
(2.5)-(2.6)  using (2.4). \hfill\carre

\bigskip

{\bf Proposition 2.2.} For all $A$ and $B$ measurable bounded functions
on $(\R^2)^I$, the function $C_h^{reg , I} (A , B)$ is bounded on
 $(\R^2)^I$ and satisfies:
 $$ \Vert C_h^{reg , I} (A , B) \Vert _{\infty} \leq \Vert A \Vert
 _{\infty} \ \Vert B \Vert _{\infty} \leqno (2.12) $$

 \bigskip

 {\it Proof.} It suffices to give the proof for $I= \{1 ,
 ..., n\}$. In view of (2.8) -(2.10), it is clear that:
$$| C_h^{reg} (A , B )(X)| \leq \Vert Op_h^{AW}(A)\Vert_{{\cal
L}(L^2(\R^n))}\ \Vert Op_h^{AW}(B)\Vert_{{\cal L}(L^2(\R^n))}$$
We also know that $ \Vert Op_h^{AW}(A)\Vert_{{\cal L}(L^2(\R^n))} \leq
\Vert A \Vert_{\infty}$. We then deduce (2.12). \hfill\carre

\bigskip

{\bf Proposition 2.3.} The following inequality holds:
$$ \Vert C_h^{hyb, I} (A , B)\Vert _{\infty } \leq  ( \pi h)^{-2|I|}  N_{I, h}(A) N_{I, h}(B)  $$
where
$$ N_{I, h}(A) =\int _{(\R^2)^I} \Vert A_{Y_I} (.) \Vert _{\infty}
dY_I$$

\bigskip

{\it Proof. }  This proposition directly follows
from (2.4) and Proposition 2.2. \hfill\carre

\bigskip

{\lessections 3. Decomposition formula.}

\bigskip

For all subsets $I$ of $\{ 1, ... , n\}$, set:
$$ T_{I, h} = \prod _{j\in I} (1 - e^{{h\over 4} \Delta_j}) \leqno (3.1)$$
$$ \Delta _I = \sum _{j \in I} \Delta_j \hskip 2cm
 \Delta_j = \partial _{x_j}^2 + \partial _{\xi_j}^2 $$
For every finite subset  $E$ of $\{1 , ..., n\}$,  ${\cal
P}_3(E)$ denotes the set of partitions of $I$ into three disjoint subsets. More precisely, an element $(I, J, K)$ of
${\cal P}_3(E)$ is an ordered sequence of three disjoint subsets
of $E$ constituting a partition of $E$, one of them or two of them possibly being empty, or even the three of them if  $E$
is itself empty.

\bigskip

{\bf Proposition 3.1.} {\it For all $A$ in $S(M,
\varepsilon , \varepsilon)$ and $B$ in $S(M', \varepsilon ,
\varepsilon )$, we have the following expression
$$ C_h^{weyl}(A,  B) = \sum _{E \subset \{1 , ... , n\} }  \sum _{(I, J
, L)\in {\cal P}_3(E)}  e^{{h\over 4} (\Delta _J+ \Delta_L) } T_{I,
h}  C_h^{hyb , E}(e^{{h\over 4} \Delta _L}  T_{J, h}A, T_{L, h}
B))\leqno (3.2)
$$

 }

 \bigskip

{\it Proof.}  We see, similarly to the paper concerning norms ([A-J-N-2]):
$$ C_h^{weyl}(A ,  B) = \sum _{I \subset \{ 1 , ... , n\} } e^{{h\over 4} \Delta_{I^c}}
T_{I, h}C_h ^{weyl}(A, B)
$$
For each subset $I$ of $\{ 1 , ... , n\} $, we also get replacing $\{ 1 , ... , n\} $ by $I^c$:
$$A = \sum _{J\subset I^c }e^{{h\over 4} \Delta_{(I \cup J)^c}}T_{J, h}
A$$
For every finite subsets $I$ and $J$, we similarly see replacing $\{ 1 , ... , n\} $ by $(I \cup J)^c$ that:
$$ B= \sum _{L\subset (I \cup J)^c }e^{{h\over 4} \Delta_{(I \cup
J\cup L)^c}}T_{L, h}B$$
Combining these three equalities yields:
$$ C_h^{weyl}(A,  B) = \sum _{E \subset \{1 , ... , n\} }  \sum _{(I, J
, L)\in {\cal P}_3(E)} e^{{h\over 4} \Delta_{I^c}} T_{I,
h}C_h^{weyl}  (e^{{h\over 4} \Delta_{(I \cup J)^c}}T_{J, h} A
,e^{{h\over 4} \Delta_{(I \cup J\cup L)^c}}T_{L, h}B)$$
Applying  definition (2.2) to each term, where $I$ is replaced by $E= I \cup J \cup L$, we obtain (3.2). \hfill\carre

\bigskip

{\lessections 4. Proof of Theorem 1.2.}

\bigskip

We first begin with the proof of  (1.3) when $\rho_j = \delta_j $ for
all $j$. These common values are denoted by $\varepsilon_j$.

\bigskip

For every subset $I$ of $\{1 , ..., n\}$,  ${\cal
M}_m(I)$ denotes the set of multi-indices  $(\alpha , \beta)$ such that:
$$ \alpha_j = \beta_j = 0 \ \ \ \ \ {\rm if} \ \ \ j\notin I \hskip
2cm \alpha_j \leq m \ \ \ \ \beta_j \leq m \ \ \ \ \ {\rm if} \ \ \
j\in I$$
For all functions $F$ in $S(M, \varepsilon ,  \varepsilon )$,
set:
$$N_{I  , h}^{(m)}(F)  = \sum _{(\alpha , \beta)\in {\cal M}_m(I)}
h^{(|\alpha|+|\beta|)/2}  \left \Vert \partial _x^{\alpha} \partial
_{\xi}^{\beta}  F \right \Vert _{L^{\infty }(\R^{2n})}  \leqno
(4.1)$$

\bigskip

{\bf Lemma 4.1.} {\it  There exists $K_0>0$ such that, for all $A$ in
$S_4(M, \varepsilon , \varepsilon )$ and $B$ in $S_4(M' ,
\varepsilon , \varepsilon )$, for every finite subset  $E$ of
$\{1 , ... , n\}$, for all $X$ in $\R^{2n} $:
$$\left  |C_h^{hyb , E }(A,  B)(X)\right |\leq  K_0^{|E|} N_{E  ,
h}^{(4)}(A) \  N_{E  , h}^{(4)}(B)\leqno (4.2)$$

}

{\it Proof.} We start from the expressions (2.5)-(2.6). Notice that, for any $j\in E$:
$$ L_jK_{E, h} (Y , Z)  = L'_jK_{E, h} (Y , Z) = K_{E, h} (Y , Z) $$
$$ L_j = \left ( 1 + {1\over h} (y_j^2  + \eta_j^2
) \right )^{-1} \ \left ( I - {h\over 4}(  \partial _{z_j}^2 +
\partial _{\zeta_j}^2 )  \right )
$$
$$ L'_j = \left ( 1 + {1\over h} (z_j^2  + \zeta_j^2
) \right )^{-1} \ \left ( I - {h\over 4}(  \partial _{y_j}^2 +
\partial _{\eta_j}^2 )  \right )
$$
This provides:
$$ C_h^{hyb, E} (A , B) (X) =    \int _{\R^{4n} }K_{E, h} (Y , Z)
\left (  \prod _{j\in E}  ( ^t L_j)^2( ^t L'_j)^2 A(X+Y)
B(X+Z)\right )
 dY dZ$$
We may write:
$$\left (  \prod _{j\in E}  ( ^t L_j)^2( ^t L'_j)^2 A(X+Y)
B(X+Z)\right ) = \sum _{(\alpha , \beta) \in {\cal M}_4(E)\atop
(\gamma , \delta ) \in {\cal M}_4(E)} \ \Phi_{\alpha , \beta ,
\gamma , \delta } \left ( {Y_I \over \sqrt h}\right ) \ \Psi_{\alpha
, \beta , \gamma , \delta } \left ( {Z_I \over \sqrt h}\right )...
$$
$$ ...  h^{(|\alpha| + |\beta |+ |\gamma |+|\delta |) /2} \Big (
\partial _{y}^{\alpha}\partial _{\eta }^{\beta}A \Big )(X+ Y)  \ \Big ( \partial
_{z}^{\gamma}\partial _{\zeta}^{\delta } B \Big ) (X+Z) $$
where all the  $\Phi_{\alpha , \beta , \gamma , \delta }$ and $\Psi_{\alpha
, \beta , \gamma , \delta }$ are functions on $(\R^2)^I$
satisfying for some universal constant  $K>0$:
$$ \int_{(\R^2)^E} | \Phi_{\alpha , \beta , \gamma , \delta
} (Y_E)| dY_E \leq K^{|E|} $$
and likewise for the $\Psi_{\alpha , \beta , \gamma , \delta }$. In particular, for $X$ fixed:
$$ C_h^{hyb , E }(A,  B)(X)  =\sum _{(\alpha , \beta) \in {\cal M}_4(E)\atop
(\gamma , \delta ) \in {\cal M}_4(E)} h^{(|\alpha| + |\beta |+
|\gamma |+|\delta |) /2}  C_h^{hyb , E } ( A _{\alpha \beta \gamma
\delta, X } ,B _{\alpha \beta \gamma \delta, X })(X)$$
with
$$A _{\alpha \beta \gamma \delta , X} =
\Phi_{\alpha , \beta , \gamma , \delta } \left ( {Y_I- X_I  \over
\sqrt h}\right ) \partial_y^{\alpha}
\partial_{ \eta}^{\beta}A$$
and similarly:
$$B _{\alpha \beta \gamma \delta , X} =
\Psi_{\alpha , \beta , \gamma , \delta } \left ( {Z_I -X_I \over
\sqrt h}\right ) \partial_z^{\gamma}
\partial_{ \zeta}^{\delta }B$$
Taking Proposition 2.3 into account, this implies
$$ \Vert  C_h^{hyb , E } (
A _{\alpha \beta \gamma \delta , X} ,B _{\alpha \beta \gamma \delta
, X }) \Vert _{\infty } \leq ( \pi h)^{-2|E|}  N_{E, h}(A_{\alpha
\beta \gamma \delta , X }) N_{E, h}(B_{\alpha \beta \gamma \delta ,
X})
$$
Besides,
$$  N_{E, h}(A_{\alpha \beta \gamma \delta , X}) \leq  \Vert \partial_y^{\alpha}
\partial_{ \eta}^{\beta}A \Vert _{\infty} \int_{(\R^2)^E} \left |
\Phi_{\alpha , \beta , \gamma , \delta } \left ( {Y_I \over \sqrt
h}\right ) \right | dY_I $$
and similarly for $B_{\alpha \beta \gamma \delta , X}$. Therefore the proof of (4.2) is completed.\hfill\carre

\bigskip

{\bf Proposition 4.2.} {\it There exists a universal constant
$K>0$ such that, for any $A$ in $S_6(M, \varepsilon ,
\varepsilon )$ and $B$ in $S_6(M', \varepsilon , \varepsilon )$,
for all $E \subset \{1 , ... , n\} $, for all $(I, J , L)\in
{\cal P}_3(E) $,
$$ \Vert e^{{h\over 4} (\Delta _J+ \Delta_L) } T_{I,
h}  C_h^{hyb , E}(e^{{h\over 4} \Delta _L}  T_{J, h}A, T_{L, h} B))
\Vert _{\infty} \leq M M' (K h) ^{|E|} \prod _{j\in E}
\varepsilon_j^2 \leqno (4.3)$$

}

\bigskip

{\it Proof.} From the definition (3.1) of $T_{I, h}$
and heat kernel properties, we get:
$$ \Vert e^{{h\over 4} (\Delta _J+ \Delta_L) } T_{I,
h}  C_h^{hyb , E}(e^{{h\over 4} \Delta _L}  T_{J, h}A, T_{L, h} B))
\Vert _{\infty} \leq (h/4)^{|I|} \left  \Vert \left [\prod _{i\in
I}\Delta_i \right ] C_h^{hyb , E}(e^{{h\over 4} \Delta _L}  T_{J,
h}A, T_{L, h} B))\right  \Vert _{\infty}$$
Using (2.5),  it is clear that, for every $i\in I$, for all $F$ and $G$:
$$ \Delta_iC_h^{hyb , E}( F , G) =C_h^{hyb , E}( \Delta_i F , G) + 2
 C_h^{hyb , E}( \partial _{x_i} F ,  \partial _{x_i} G)+ 2
 C_h^{hyb , E}( \partial _{\xi_i} F ,  \partial _{\xi _i} G)
+ C_h^{hyb , E}(  F , \Delta_i G)$$
This yields:
$$  \left [\prod _{i\in
I}\Delta_i \right ] C_h^{hyb , E}(  F, G) = \sum _{(\lambda , \mu ,
\lambda ', \mu')\in \widetilde  {\cal M}(I)} c_{\lambda , \mu ,
\lambda ', \mu'} C_h^{hyb , E}(
\partial _x^{\lambda} \partial _{\xi}^{\mu}  F , \partial _x^{\lambda '}
\partial _{\xi}^{\mu '}G) \leqno (4.4)$$
where $\widetilde {\cal M}(I)$ denotes the set of multi-indices
$(\lambda , \mu , \lambda ', \mu')$ such that:
$$ \lambda _i = \mu_i = \lambda '_i = \mu'_i = 0 \ \ \ \ {\rm if} \
\ \ \ \ i \notin I \hskip 2cm \lambda _i + \mu_i + \lambda '_i +
\mu'_i = 2  \ \ \ \ {\rm if} \ \ \ \ \ i \in I $$
In (4.4), the $c_{\lambda , \mu , \lambda ', \mu'}$ are constants with absolute values smaller or equal than
 $2 ^{|I|}$. According to Lemma 4.1 we deduce that:
$$\Vert e^{{h\over 4} (\Delta _J+ \Delta_L) } T_{I,
h}  C_h^{hyb , E}(e^{{h\over 4} \Delta _L}  T_{J, h}A, T_{L, h} B))
\Vert _{\infty} \leq ... $$
$$ ... \leq 2 ^{|I|}  K_0^{|E|}  \sum _{(\lambda , \mu , \lambda ', \mu')\in \widetilde {\cal
M}(I)}  N_{E  , h}^{(4)}(\partial _x^{\lambda}
\partial _{\xi}^{\mu} e^{{h\over 4} \Delta _L}  T_{J, h} A) \
N_{E  , h}^{(4)}(\partial _x^{\lambda '}
\partial _{\xi}^{\mu '} T_{L, h} B) $$
Besides, for all $A$  in $S_6 (M, \varepsilon , \varepsilon )$,
for every disjoint subsets  $I$, $J$ and $L$ of $E\subset \{1
, ..., n\}$, for each multi-index $(\lambda , \mu , \lambda ',
\mu')$ in $\widetilde {\cal M}(I)$ and for each multi-index
$(\alpha , \beta)$ in ${\cal M}_4(E)$, if $h \varepsilon_j^2 \leq
1$ for all $j\leq n$:
$$h^{(|\alpha | + |\beta |)/2}  \Vert  \partial_x^{\lambda }
\partial_{ \xi}^{\mu}  \partial_x^{\alpha}
\partial_{ \xi}^{\beta} e^{{h\over 4} \Delta _L}  T_{J, h}A \Vert _{\infty}
\leq M (h/4)^{|J|} \prod _{j \in J} \varepsilon_j^2 \prod _{i \in I}
\varepsilon_i^{\lambda _i + \mu_i}  $$
Similarly, if $B$ is in $S_6 (M', \varepsilon , \varepsilon )$,
for each multi-index $(\lambda , \mu , \lambda ', \mu')$ in
$\widetilde {\cal M}(I)$ and for each multi-index $(\gamma ,
\delta)$ in ${\cal M}_4(E)$
$$\Vert  \partial_x^{\lambda ' }
\partial_{ \xi}^{\mu ' }   \partial_x^{\gamma}
\partial_{ \xi}^{\delta} T_{L, h} B \Vert _{\infty}\leq
M' (h/4)^{|L|} \prod _{\ell  \in L} \varepsilon_{\ell }^2 \prod _{i
\in I} \varepsilon_i^{\lambda' _i + \mu'_i} $$
The numbers of elements of $\widetilde {\cal M} (I)$ (for $I
\subset E$) and those  of ${\cal M}_4(E)$ are both being bounded  by $K^{|E|}$. Thus, we indeed deduce (4.3). \hfill\carre

\bigskip

{\it End of the proof of Theorem 1.2.} In view of
Propositions 3.1 and 4.2, if  $A$ belongs to $S_6 (M, \varepsilon
, \varepsilon )$ and $B$  lies in $S_6 (M', \varepsilon , \varepsilon
)$ then:
$$ \Vert C_h^{weyl}(A,  B) \Vert _{\infty} \leq  \sum _{E \subset \{1 , ... , n\} }  \sum _{(I, J
, L)\in {\cal P}_3(E)}M M' (K h) ^{|E|} \prod _{j\in E}
\varepsilon_j^2
$$
For every finite subset $E$ of $\{ 1 , ... , n\}$, the number of elements of ${\cal P}_3(E)$ is $3! ( 1 + \sigma_p^2 +
\sigma_p^3)$ where $p= |E|$ and the  $\sigma_p^k$ are the Stirling numbers of second kind. This number of elements is bounded by $K_1^{|E|}$. Consequently:
$$ \Vert C_h^{weyl}(A,  B) \Vert _{\infty} \leq  M M' \sum _{E \subset \{1 , ... , n\}
}( K K_1 h)^{|E|} \prod _{j\in E} \varepsilon_j^2 $$
$$ = M M' \prod _{1\leq j \leq n} ( 1 + K K_1 h \varepsilon_j^2)$$
proving the first claim of Theorem 1.2 when
$\rho_j = \delta _j = \varepsilon_j$ for all $j$. In the general case, we set, for any function $F$ on $\R^{2n}$ and for
any sequence $\lambda = (\lambda _1 , ..., \lambda_n)$ of positive real numbers:
$$ (\delta_{\lambda }F) (x , \xi) = F \left ( \lambda _1 x_1 , ...,
\lambda _n x_n , {\xi _1 \over \lambda _1}, ... ,{\xi _n \over
\lambda _n}\right ) $$
In particular,
$$ \delta _{\lambda}C_h^{weyl} (A , B) = C_h^{weyl} ( \delta_{\lambda
} A , \delta _{\lambda }B)$$
If $A$ belongs to $S_6(M, \rho, \delta)$ and $B$ is in $S_6(M', \rho,
\delta)$, the two sequences $(\rho_j)$ and $(\delta_j)$ being positive, then $\delta_{\lambda }A$ lies in $S_6(M,
\varepsilon , \varepsilon)$ and $\delta_{\lambda }B$ in $S_6(M',
\varepsilon , \varepsilon)$, when setting $\varepsilon_j = \sqrt
{\rho_j \delta _j}$ and $\lambda_j=\sqrt { \delta_j / \rho_j}$. The preceding result applied to  $\delta_{\lambda }A$ and
$\delta_{\lambda }B$ allows to deduce a bound in the supremum norm
of $\delta _{\lambda} C_h^{weyl}(A, B)$, which is the same as the one
of $C_h^{weyl} (A , B)$. The first claim in Theorem
1.2 is therefore derived if all the $\rho_j$ and $\delta_j$ are positive and also proved by continuity if some of them are vanishing.   For the second claim in Theorem 1.2, we remark that, if $(\alpha , \beta)$ is a multi-index such that
$\alpha _j \leq m$ and $\beta _j \leq m$ for all $j$:
$$ \partial _x^{\alpha} \partial _{\xi} ^{\beta} C_h ^{weyl } (A ,
B)= \sum _{ \alpha ' + \alpha '' = \alpha \atop \beta ' + \beta '' =
\beta } { \alpha ! \over \alpha ' ! \alpha '' !} {\beta ! \over
\beta '! \beta '' !} C_h^{weyl} ( \partial _x^{\alpha '}
\partial_{\xi}^{\beta '} A , \partial _x^{\alpha ''}
\partial_{\xi}^{\beta ''} B ) $$
If $A$ is in $S_{m+6}(M, \rho, \delta)$ and $B$ in $S_{m+6}(M',
\rho, \delta)$ then  $\partial _x^{\alpha '}
\partial_{\xi}^{\beta '} A $ is in $S_6 (M \prod \rho_j ^{\alpha' _j}
\delta _j^{\beta' _j  }, \rho, \delta  )$ and similarly for the other factor. Inequality (1.3) gives:
$$\Vert  \partial _x^{\alpha} \partial _{\xi} ^{\beta} C_h ^{weyl } (A ,
B)\Vert _{\infty} \leq \left [ \sum _{ \alpha ' + \alpha '' = \alpha
\atop \beta ' + \beta '' = \beta } { \alpha ! \over \alpha ' !
\alpha '' !} {\beta !  \over \beta '! \beta '' !} \right ] M'' \prod
\rho_j ^{\alpha _j} \delta _j^{\beta _j} $$
The above sum is bounded by $2^{|S(\alpha )|+ |S(\beta)|}$
where $S(\alpha)$ is the support of $\alpha$, i.e., the set of indices
$j$ satisfying $\alpha_j \not= 0$. We therefore deduce the second claim in Theorem 1.2. \hfill\carre

\bigskip

{\lessections 5. Proof of Theorem 1.3.}

\bigskip

In this last section, the remainder term of the semiclassical expansion of the Weyl composition
is considered.

\bigskip

{\bf Proposition 5.1.} {\it Under the hypotheses of Theorem 1.3, the
function $R_N(\cdot  ,  h)$ defined in (1.4) may be written as:
$$R_N(X ,  h) ={N h^N \over (2i)^N }\sum_{ |\alpha |+|\beta | = N}
{(-1)^{|\beta|} \over \alpha ! \beta !} \int_0^1 (1- \theta)^{N-1}
 C_{\theta h}  ^{weyl} (\partial _x ^{\beta}
\partial_{\xi}^{\alpha} A , \partial _x^{\alpha}
\partial _{\xi}^{\beta} B ) (X) d\theta \leqno (5.1)$$
}

\bigskip

{\it Proof.} Set:
$$ F_h (A , B , X , \theta) =  C_{\theta h}  ^{weyl} (A , B) (X) $$
From (1.1), it is clear that:
$$ F_h (A , B , X , \theta) =   \int _{\R^{4n}} A(X+ Y) B(X+Z) K_h (Y , Z , \theta ) dY
dZ$$
$$ K_h (Y , Z , \theta ) = ( \pi \theta h)^{-2n}e^{-{2i \over
\theta  h} \sigma ( Y , Z) }$$
We then verify that:
$$ {\partial K_h \over \partial \theta} = L K_h
\hskip 2cm L = {h\over 2i} \sigma (\nabla _Y , \nabla _Z)$$
In particular:
$$\partial _{\theta}^m F_h (A , B , X , \theta)= \int _{\R^{4n}}
K_h (Y , Z , \theta )  L^m [A(X+ Y) B(X+Z)] dY dZ $$
We have $K_h(\cdot  , \cdot , 0) = \delta _{(0, 0)}$. Then, we may write:
$$ C_h ^{weyl} (A , B) (X) =  \sum _{k=0}^{N-1}
{h^k \over (2i)^k k!} \sigma (\nabla _Y , \nabla _Z) ^k [A(X+Y)
B(X+Z)]\Bigg |_{Y= Z = 0} + R_N(X ,  h)
$$
$$R_N(X ,  h) ={h^N \over (2i)^N (N-1)!}  \int _{\R^{4n}\times [0, 1]}
(1 - \theta)^{N-1} K_h (Y , Z , \theta )  \sigma (\nabla _Y , \nabla
_Z) ^N [A(X+Y) B(X+Z)]  d\lambda (Y, Z)d\theta
$$
Besides:
$$ {1 \over N!}   \sigma (\nabla _Y , \nabla
_Z) ^N [A(X+Y) B(X+Z)] = \sum_{ |\alpha |+|\beta | = N}
{(-1)^{|\beta|} \over \alpha ! \beta !} \Big [ \partial _x ^{\beta}
\partial_{\xi}^{\alpha} A(X+Y) \Big ] \ \Big [ \partial _x^{\alpha}
\partial _{\xi}^{\beta} B( X+Z) \Big ]$$
From these considerations, we then deduce (5.1). \hfill\carre

\bigskip

{\it Proof of  Theorem 1.3.} If $A$ belongs to $S_{m} (M, \rho ,
\delta )$, if  $B$ is in $S_{m} (M', \rho , \delta )$ ($m\geq 6$),
and if $|\alpha |+|\beta | = N$,  then
$$\partial _x ^{\beta} \partial_{\xi}^{\alpha} A\in S_{m-N}\left (
 M
\prod_{j\leq n}  \rho _j^{\beta_j}  \delta _j^{\alpha_j}  , \rho,
\delta \right )$$
According to Theorem 1.2, we then deduce that:
$$ C_{\theta h}  ^{weyl} (\partial _x ^{\beta}
\partial_{\xi}^{\alpha} A , \partial _x^{\alpha}
\partial _{\xi}^{\beta} B ) \in S _{m - N - 6}  \left ( M''\prod_{j\leq n}
(\rho_j \delta _j)^{\alpha_j + \beta _j} , 2 \rho , 2 \delta \right
)
$$
Noticing that
$$ N \int_0^1 (1 - \theta )^{N-1} d\theta = 1$$
we obtain:
$$R_N(\cdot  ,  h)\in S_{m - N - 6}\left ( M''{h^N \over 2^N }\sum_{ |\alpha |+|\beta | = N}
{1 \over \alpha ! \beta !}\prod_{j\leq n} (\rho_j \delta
_j)^{\alpha_j + \beta _j} , 2 \rho , 2 \delta \right ) $$
$$ =S_{m - N - 6}\left ( M'' {h^N \over N! } \left [  \sum _{j=1}^n \rho_j
\delta _j\right ]^N , 2 \rho , 2 \delta \right ) $$
and the proof of Theorem 1.3 is completed.\hfill\carre
\bigskip

\bigskip
\centerline{\lessections References}

\bigskip

[A-J-N-1] L. Amour, L. Jager, J. Nourrigat, {\it Bounded Weyl
pseudodifferential operators in Fock space,} preprint,
arXiv:1209.2852, sept. 2012.

\smallskip

[A-J-N-2] L. Amour, L. Jager, J. Nourrigat, {\it On bounded
pseudodifferential operators in a high-dimensional setting,}
preprint, arXiv:1303.1972, march 2013, submitted to Proceedings of
the A.M.S, under revision.

\smallskip

[BE] F.A. Berezin, {\it Quantization,} (Russian), Izv. Akad. Nauk.
SSSR, Ser. Mat, {\bf 38} (1974), 1116-1175.

\smallskip

[F] G. B. Folland, {\it Harmonic analysis in phase space.} Annals of
Mathematics Studies, {\bf 122}. Princeton University Press,
Princeton, NJ, 1989.

\smallskip

[HO] L. H\"ormander, {\it The analysis of linear partial
differential operators,} Volume III, Springer, 1985.

\smallskip

[LA-1] B. Lascar, {\it Noyaux d'une classe d'op\'erateurs
pseudo-diff\'erentiels sur l'espace de Fock, et applications.}
S\'eminaire Paul Kr\'ee, 3e ann\'ee (1976-77), Equations aux
d\'eriv\'ees partielles en dimension infinie, Exp. No. 6, 43 pp.

\smallskip

[LA-2] B. Lascar,  {\it Une classe d'op\'erateurs elliptiques du
second ordre sur un espace de Hilbert,}  J. Funct. Anal. {\bf 35}
(1980), no. 3, 316-343.

\smallskip

[LA-3] B.  Lascar, {\it  Op\'erateurs pseudo-diff\'erentiels en
dimension infinie. Applications. } C. R. Acad. Sci. Paris S\'er. A-B
{\bf 284} (1977), no. 13, A767-A769,

\smallskip

[LA-4] B. Lascar, {\it  Op\'erateurs pseudo-diff\'erentiels d'une
infinit\'e de variables, d'apr\`es M. I. Visik.} S\'eminaire Pierre
Lelong (Analyse), Ann\'ee 1973-1974, pp. 83–90. Lecture Notes in
Math., {\bf 474}, Springer, Berlin, 1975.

\smallskip

[LER] N. Lerner, {\it Metrics on the phase space and non
self-adjoint pseudo-differential operators,} Birkh\"auser Springer,
2010.

\smallskip

[M] A. Martinez, {\it An introduction to semiclassical and
microlocal analysis.}  Universitext, Springer-Verlag, New York,
2002.

\smallskip

[R] D. Robert, {\it Autour de l'approximation semi-classique},
Progress in Mathematics 68, Birkh\"auser Boston, Inc., Boston, MA,
1987.

\smallskip

[S] M.A. Shubin, {\it Pseudodifferential operators and spectral
theory,} translated from the 1978 Russian original, Springer, 2001.

\smallskip

[U-1] A. Unterberger, {\it  Oscillateur harmonique et op\'erateurs
pseudo-diff\'erentiels},   Ann. Inst. Fourier (Grenoble) {\bf 29}
(1979), no. 3, xi, 201-221.

\smallskip

[U-2] A. Unterberger, {\it Les op\'erateurs m\'etadiff\'erentiels},
in Complex analysis, microlocal calculus and relativistic quantum
theory, Lecture Notes in Physics {\bf 126} (1980) 205-241.

\bigskip

{\it Address:} LMR EA 4535 et FR CNRS 3399, Universit\'e de Reims
Champagne-Ardenne, Moulin de la Housse, B. P. 1039, F-51687
Reims, France.

\medskip

{\it Email:} laurent.amour@univ-reims.fr
\smallskip
{\it Email:}  jean.nourrigat@univ-reims.fr

\end